\documentclass[12pt,reqno]{amsart}
\usepackage[latin1]{inputenc}  
\usepackage{amssymb}
\usepackage{stmaryrd}
\usepackage[all]{xy}
\usepackage{graphicx}
\usepackage{hyperref}
\theoremstyle{plain}

\newtheorem{theo}{Theorem}
\newtheorem{lemme} {Lemma}

\newtheorem{prop}{Proposition}
\newtheorem{coro}{Corollary}

\newcommand{\finpreuve}{\mbox{} \hfill \mbox{$\Box$}}
\newcommand{\finpreuvelemme}{\mbox{} \hfill \mbox{$\blacksquare$}}
\newenvironment{preuve}{\noindent {\it Proof :}}{\finpreuve}
\newenvironment{preuvelemme}{\noindent {\it Proof of the lemma:
}}{\finpreuvelemme}

\makeatletter
\newcounter{enum@ux}
\def\brkenum#1{%
\setcounter{enum@ux}{\value{enum\romannumeral\the\@enumdepth}}%
\end{enumerate} #1%
\begin{enumerate}%
\addtocounter{enum\romannumeral\the\@enumdepth}{\value{enum@ux}}}
\makeatother

\makeindex

\newcommand{\pg}{\mathfrak{p}}

\title{Generalised Mertens and Brauer--Siegel Theorems}
\author{Philippe Lebacque}
\date{}
\address{Institut de Math\'ematiques de Luminy, Université de la M\'editerran\'ee, Campus de Luminy, case 907, 13288 Marseille Cedex 9, France}
\email{lebacque@iml.univ-mrs.fr}
\begin{document}
\maketitle

\section{Introduction}
In this article, we prove a generalisation of the Mertens theorem for prime numbers to number fields and algebraic varieties over finite fields, paying attention to the genus of the field (or the Betti numbers of the variety), in order to make it tend to infinity and thus to point out the link between it and the famous Brauer--Siegel theorem. Using this we deduce an explicit version of the generalised Brauer--Siegel theorem under GRH, and a unified proof of this theorem for asymptotically exact families of almost normal number fields.

The classical Brauer--Siegel theorem is a well-known theorem which describes the behaviour of the quantity $hR$ (the product of the class number and the regulator) in a family of number fields with growing genus under the conditions that the genus grows much faster than the degree and some additional properties like normality or the Generalised Riemann Hypothesis (GRH) to deal with the Siegel zeroes. These two hypotheses are of different nature: omitting the first one changes the final result, while the second one is a technical hypothesis. Tsfasman and Vl\u adu\c t \cite{TVF} were able to remove the first hypothesis, which led to the so called generalised Brauer--Siegel theorem, and Zykin \cite{Zyk} was able to replace "normality" by "almost normality" in the second one using results of Stark and Louboutin. He also managed to generalise the Brauer--Siegel theorem to the case of smooth absolutely irreducible projective varieties over finite fields.

As for the Mertens theorem, proven by Mertens in the case of $\mathbb{Q}$, and much later generalised by Rosen \cite{ROS} both in cases of number and function fields, it can be regarded as the Brauer--Siegel theorem in the finite steps of the family. An explicit Mertens theorem leads therefore to an explicit formulation of the generalised Brauer--Siegel theorem. We first recall the formulations of the (generalised) Brauer--Siegel theorem and Mertens theorem, then we prove their explicit versions for number fields and smooth projective absolutely irreducible varieties over finite fields, and finally we deduce the explicit generalised Brauer Siegel theorem.

\section{Around the Brauer--Siegel theorem}

Let us now recall the notations and the definitions involved in the generalized Brauer--Siegel theorem, and state it for global fields and smooth absolutely irreducible projective algebraic varieties ($s.a.i.p.a.v.$) over the finite field $\mathbb{F}_r.$
Throughout this paper we will write $(NF)$ and $(V)$ to say that something is true in the case of number fields and $s.a.i.p.a.v.$ respectively.

\subsection{Number field case}

Given a number field $K$, let $\zeta_K$ be the usual zeta-function of the field $K$ and $\varkappa_K$ be its residue at $s=1.$ Denote by $\Phi_q(K)$ the number of places of $K$ whose norm is equal to $q.$ Let $(K_i)_{i\in\mathbb{N}}$ be a family of finite extensions of $\mathbb{Q}$. Denote by $g_i:=\log{\sqrt{|\text{Discr}(K_i)}|}$ the genus of $K_i$ and $n_i$ its degree. Recall that $(K_i)_{i\in\mathbb{N}}$ is said to be asymptotically exact if $\phi_q:=\lim\Phi_q(K_i)/g_i$ exist for all prime powers $q$ and if $\phi_\mathbb{R}:=\lim r_1(K_i)/g_i$ and $\phi_\mathbb{C}:=\lim r_2(K_i)/g_i$ exist, where $r_1(K_i)$ and $r_2(K_i)$ stand for the number of real and complex places of $K_i$ respectively. We put $\phi_\infty=\phi_\mathbb{R}+2\phi_\mathbb{C}.$ Being asymptotically exact is not a restrictive property. In fact, every tower of global fields is asymptotically exact, and each family of number fields contains an asymptotically exact subfamily. In the classical Brauer--Siegel theorem, all these $\phi_q$ are equal to zero because of the assumption $n_i/g_i\to 0:$

\begin{theo}[Classical Brauer--Siegel]\label{CBS}
Assume that the family of number fields $(K_i)$ is normal over $\mathbb{Q}$ or that \emph{GRH} holds, and assume that $\lim_i n_i/g_i=0$. Then $\log{h_iR_i}\sim g_i.$
\end{theo}

In order to prove this theorem, we need to use the class number formula: $$\varkappa_K=\frac{2^{r_1}(2\pi)^{r_2}}{w|d_K|^{1/2}}hR.$$ The result can be reformulated in this way: $$\lim_i\frac{\log\varkappa_{K_i}}{g_i}=0.$$ Suppressing this hypothesis leads to the Tsfasman-Vl\u adu\c t Brauer--Siegel theorem (T-V B-S). This time, the $\phi_q$ are not always equal to zero: 

\begin{theo}[T-V BS (2002)] \label{GBS}Let $(K_i)_{i\in\mathbb{N}}$ be an asymptotically exact family of number fields. Assume either that \emph{GRH} holds, or that $(K_i)$ is a family of almost-normal number fields (We will say that a number field $L$ is almost normal if there exists a tower $L_0\subset\dots\subset L_n=L$ of fields such that $L_{i+1}$ is normal over $L_i$ for all $i$). Then the limit $\kappa=\lim_i \frac{\log{\varkappa_{K_i}}}{g_i}$ exists and we have the following equality: $$\kappa=\sum_{q} \phi_q \log\left(\frac{q}{q-1}\right)<+\infty,$$ where the sum is taken over all the powers of prime numbers. 
\end{theo}

In their paper \cite{TVF}, Tsfasman and Vl\u adu\c t proved this theorem without the assumption of GRH for asymptotically good families of almost normal number fields (this means $\lim n_i/g_i>0$), and Zykin \cite{Zyk} proved this is also true for asymptotically bad families. In order to get this result, we have to deal with two inequalities, but one of them is always satisfied :

\begin{theo}[BS Inequality]\label{IBS}
 Let $\mathcal{K}=(K_i)_{i\in\mathbb{N}}$ be an asymptotically exact family of number fields. Then $$\limsup_i \frac{\log{\varkappa_{K_i}}}{g_i}\leq\sum_q \phi_q \log\left(\frac{q}{q-1}\right)<+\infty.$$
\end{theo}

The difficulties come from the second inequality 
$$\sum_q \phi_q \log\left(\frac{q}{q-1}\right)\leq\liminf_i \frac{\log{\varkappa_{K_i}}}{g_i},$$ which requires technical assumptions.

\subsection{Case of algebraic varieties over a finite field}

Consider an algebraic variety $X$ of dimension $d$, defined over a finite field $\mathbb{F}_r.$ Suppose that $X$ is smooth, projective and absolutely irreducible and let $|X|$ denote the set of its closed points. For $\mathfrak{p}\in|X|$ and $k(\pg)$ its residue field, let $\deg(\pg)$ be the degree of the field extension $[k(\pg):\mathbb{F}_r].$ Define now for $m\geq 1$ the $\Phi$-numbers as before: $$\Phi_{r^m}:=\#\left\{\mathfrak{p}\in|X|\  |\ \text{deg}(\mathfrak{p})=m  \right\}.$$
Put $\bar{X}=X\otimes\mathbb{F}$ where $\mathbb{F}=\bar{\mathbb{F}_r}$ is the algebraic closure of $\mathbb{F}_r.$ Let $\ell$ be a prime different from $p$. Let $b_i=\text{dim}H^i(\bar{X},\mathbb{Q}_\ell)$ for $0\leq i\leq 2d$ be the Betti numbers for the $\ell$-adic etale cohomology of $X$. $X$ is smooth, so they do not depend on $\ell$ and verify the equality $b_i=b_{2d-i}$ because of the Poncar\'e duality.  Let $b_X=\max_{i=0,..,2d}b_{i}.$ In the case of dimension $1,$ $b_0=b_2=1$ and $b_1=g,$ so we have $b_X=\max(g,1).$ In this theory the quantity $b_X$ will play the role of the genus of number fields (and function fields). Since the asymptotic theory of varieties of dimension higher than $1$ is not yet well understood, we do not know exactly which quantity is the exact analogue of the genus. We chose this number $b_X$ because it was easier to compute the sums, but it might happen that the sum of $b_i$ or a certain sum of $b_i$'s with coefficients depending on $r$ could make a better choice. However, unless we want to make $r$ or $d$ grow, all these choices are equivalent. 

By the famous Deligne--Grothendieck theorem, the zeta function of $X$ verifies $$Z(X,t)=\prod_{i=0}^{2d}P_i(t)^{(-1)^{i+1}},$$ where $$P_i(t)=\prod_{j=1}^{b_i}(1-\omega_{i,j}r^{i/2}t),$$
$\omega_{i,j}$ being algebraic numbers of module $1$ and $P_0(t)=1-t,$ $P_{2d}(t)=1-r^dt$.
 We will consider $\zeta_X(s)=Z(X,r^{-s}),$ and $\varkappa_X=\rm{Res}\it _{s=d}\zeta_X.$
 
 Let us fix the dimension $d,$ and let $X$ go through a family of  $s.a.i.p.a.v.$ of dimension $d.$ We say that the family $(X_i)_{i\in\mathbb{N}}$ is asymptotically exact if $b_{X_i}\to\infty$, and if, for all $m\geq 1,$ the limit $\phi_{r^m}=\lim_{i}\Phi_{r^m}/b_{X_i}$ exist.
 
 We can now formulate a generalisation of the Brauer--Siegel theorem for varieties of dimension $d.$ It was proved by Tsfasman and Vl\u adu\c t in the function field case \cite{TVF}, and by Zykin (unpublished) in the case of $d>1$, using a different definition of $b_X.$
 \begin{theo} \label{VBS} Let $(X_i)_{i\in\mathbb{N}}$ be an asymptotically exact family of s.a.i.p.a.v. of dimension $d$ defined over $\mathbb{F}_r.$ Then $\kappa=\lim_i\log(\varkappa_{X_i})/b_{X_i}$ exists and we have the following equality: $$\kappa=\sum_{m=1}^\infty\phi_{r^m}\log\left(\frac{r^{dm}}{r^{dm}-1}\right)$$
 \end{theo}

 Unfortunetaly, we do not know any reasonable interpretation of the residue of the zeta function at $s=d,$ such as we have for $s=1$ through the class number formula in the number field and function field cases .

\section{Mertens theorem and its relation to the generalised BS theorem}

If one wants to get an explicit version of the generalized Brauer--Siegel equality, one need to know what happens explicitly between the residue $\varkappa_{K_i}$ and $\sum_{q\leq x}\phi_q\log\frac{q}{q-1}$ at the finite steps of the family. This is given by the Mertens theorem. 

\begin{theo}[Mertens]  

\begin{align*}
(V)\qquad \qquad &\ \, \prod_{\stackrel{P\in |X|}{\deg(P)\leq N}}\left(1-\frac{1}{\mathcal{N}P^d}\right)=\frac{e^{-\gamma_{X}}}{N}+\mathcal{O_X}\left(\frac{1}{N^2}\right),\\
(NF) \qquad \qquad & \ \ \, \prod_{\stackrel{P\in P_f(K)}{\mathcal{N}P\leq x}}\left(1-\frac{1}{\mathcal{N}P}\right)=\frac{e^{-\gamma_{K}}}{\log x}+\mathcal{O_K}\left(\frac{1}{\log^2{x}}\right),\\
(NF-GRH) \qquad \qquad & \ \ \, \prod_{\stackrel{P\in P_f(K)}{\mathcal{N}P\leq x}}\left(1-\frac{1}{\mathcal{N}P}\right)=\frac{e^{-\gamma_{K}}}{\log x}+\mathcal{O_K}\left(\frac{1}{\sqrt{x}}\right),\\
\end{align*}

where \begin{align*}
(V) \qquad\qquad & \gamma_{X}=\gamma+\log{(\varkappa_X\log{r})},\\
(NF)\qquad\qquad  & \gamma_{K}=\gamma+\log{\varkappa_K},
\end{align*}
$P_f(K)$ being the set of finite places of $K,$ and $\mathcal{N}P$ denoting the absolute norm of the place $P.$
\end{theo}

The function field case and the number field case are due to Rosen, who proved them following the classical proof of the Mertens classical theorem \cite{HW}. But he paid no attention to the behaviour of the constants in field extensions. Unfortunately we did not know about his work before having ended ours. Mireille Car also proposed in \cite{Car} a different proof in the case of function fields. In the number field case, we also follow the classical Mertens proof with small variations in order to get an explicit version of this theorem, which takes into account the genus and the degree of $K$. In the case of varieties over finite fields, we present a natural proof using explicit formulae. We prove in fact the following sharper results:

Without assuming GRH, we have to deal with exceptional zeroes. A real zero $\rho$ of $\zeta_K$ is said to be exceptional if $1-(8\,g)^{-1}\leq\rho<1.$ A number field has at most one exceptional zero. A real zero $\rho$ is a Siegel zero if $1-(32\,g)^{-1}\leq\rho<1.$
\begin{theo}Let $K$ be a number field. $$(NF)\qquad \sum_{q\leq x}\Phi_q\log\left(\frac{q}{q-1}\right)=\log\log{x}+\gamma+\log{\varkappa_K}+\tau_1(x)+\frac{1}{1-\rho}\tau_2(x), $$ and there exist effective constants $C,c_1,c_2$ such that, for all $x\geq Cn\,g^2,$ $$|\tau_1(x)|\leq c_1\frac{1}{\log{x}},$$ 
\begin{align*}|\tau_2(x)|&\leq \displaystyle c_2\frac{1}{\log{x}}, \text{ if K has an exceptional zero } \rho,\\ 
 &=0 \text{ otherwise.}\end{align*}
\end{theo}
The condition on $x$ does not allow us to have explicit results as in the case where GRH holds, but these results, combined with theorem \ref{IBS}, lead us to the unified proof of the Brauer--Siegel theorem, and to other nice results around the Brauer--Siegel theorem and the family of $\phi_q.$

\begin{coro} Let $(K_i)$ be an asymptotically exact family of almost normal number fields. Then the limit $\lim_{i\to\infty}\frac{\log{\varkappa_{K_i}}}{g_{K_i}}=\kappa$ exists and verifies the following equality:
$$\sum_{q}\phi_q\log{\frac{q}{q-1}}=\kappa,$$  the sum being taken over all prime powers $q$.
\end{coro}
We cannot suppress the hypothesis of normality, because of exceptional zeroes that can appear in the family. But we can say something more in the general case:

\begin{prop}\label{propcomp}Let $(K_i)$ be an asymptotically exact family of number fields.
\begin{enumerate}\item Assume that $\lim_i \frac{n_i\log{n_i}}{g_i}=0$. Then $\kappa$ exists and equals $0.$ 
\item Assume that the family $(K_i)$ is asymptotically good ( i.e. $\phi_\infty>0$ ), and that there are infinitely many Siegel zeroes in the family. Then $$\sum_q\phi_q\log\frac{q}{q-1}\leq\phi_\infty\log{\frac{e}{\phi_\infty}}.$$
\end{enumerate}

\end{prop}

If we assume GRH, there is no condition on $x,$ and we have the following result, which leads to an explicit version of the generalised Brauer--Siegel theorem.
\begin{theo} [GRH Mertens theorem] Assume that \emph{GRH} holds. Then 
 \begin{align*}
 (V)\, \ \qquad\ \qquad \sum_{m=1}^{N}\Phi_{r^m}\log\left(\frac{r^{dm}}{r^{dm}-1}\right)=&\log N+\gamma+\log{(\varkappa_X\log{r})}\\
 &+\mathcal{O}\left(\frac{1}{N}\right)+b_X\,\mathcal{O}\left(\frac{r^\frac{-N}{2}}{N}\right),\\
(NF-GRH) \qquad \sum_{q\leq x}\Phi_q\log\left(\frac{q}{q-1}\right)=&\log\log{x}+\gamma+\log{\varkappa_K}\\
&+n_K\,\mathcal{O}\left(\frac{\log{x}}{\sqrt{x}}\right)+g_K\,\mathcal{O}\left(\frac{1}{\sqrt{x}}\right),
\end{align*}
where the constants involved in the $\mathcal{O}$ are effective and do not depend either on $X,$ or on $K.$
\end{theo}

\begin{coro}\label{EBS}Let $(K_i)$ be an asymptotically exact family of number fields, and $(X_i)$ of s.a.p.a.i.v. of dimension $d.$ Assume \emph{GRH} in the number field case. Then the limit $\lim_{i\to\infty}\frac{\log{\varkappa_{K_i}}}{g_{K_i}}=\kappa$ (resp. $\lim_{i\to\infty}\frac{\log{\varkappa_{X_i}}}{g_{X_i}}=\kappa$) exists, and we have:  
\begin{align*}
(V)\qquad \qquad \sum_{q\leq r^N}\phi_q\log{\frac{q}{q-1}}=& \kappa+\mathcal{O}\left(\frac{r^\frac{-N}{2}}{N}\right),\\
(NF-GRH)\qquad \qquad \sum_{q\leq x}\phi_q\log{\frac{q}{q-1}}=& \kappa+\mathcal{O}\left(\frac{\log{x}}{\sqrt{x}}\right).
\end{align*}
\end{coro}

\section{Proof of the Mertens theorem}

\subsection{Proof in the number field case}

In order to prove the Mertens theorem for number fields, we follow the nice proof of the classical Mertens theorem of \cite{HW} as Rosen does in his article, but we use another counting function for prime ideals. In addition, we need a precise version of the Mertens theorem, so we will have to do the work once again, sketching Rosen's proofs.

Let $K$ be a number field, let $n=[K:\mathbb{Q}]$ be its degree and $g=\frac{1}{2}\log|{\text{Discr}(K)}|$. Put $\pi(x):=\#\left\{P\in P(K)\,|\,\mathcal{N}P\leq x\right\}$. One can estimate $\pi(x)$ by the following bound due to Lagarias and Odlyzko, and improved by Serre \cite{SC}:

Consider the Li-function defined by $$\rm Li\it (x)=\int_2^x \frac{dt}{\log t}.$$
\begin{theo}[Prime ideals theorem] 
$$(NF) \qquad \qquad  \pi(x)=\rm Li\it(x)+\Delta(x),$$ where, for all $x$ such that $$ (C1) \qquad\ \log{x}\geq c_3\, n\, g^2 $$

$$|\Delta(x)|\leq \rm Li\it(x^\rho)+c_1 x \exp{\left(-c_2\,n^{-\frac{1}{2}}\log^\frac{1}{2}(x)\right)},$$the term in $\rm Li\it(x^\rho)$ is only there if $\zeta_K$ has an exceptional zero $\rho$. Under \emph{GRH}, one has the stronger result available for all $x\geq 2$:
$$(NF-GRH)\qquad \qquad   |\Delta(x)|\leq c\, x^\frac{1}{2}\left(2g+n\log{x}\right).$$
\end{theo}


First, we will give an asymptotic expression for $\sum\frac{1}{\mathcal{N}P}:$

\begin{prop}
$$\sum_{\mathcal{N}P\leq x} \frac{1}{\mathcal{N}P}=\log\log{x}+B+o(1).$$
\end{prop}

\begin{preuve} 
We have the formula: $$C(x)=\sum_{\mathcal{N}P\leq x} \frac{1}{\mathcal{N}P}=C(2)+\int_2^{x}\frac{d \pi(t)}{t}=C(2)+\int_2^{x}\frac{dt}{t\log{t}}+\int_2^x \frac{d\Delta(t)}{t}, $$
thus $$\sum_{\mathcal{N}P\leq x} \frac{1}{\mathcal{N}P}=\int_2^{x}\frac{dt}{t\log{t}}+\frac{\Delta(x)}{x}+\int_2^x \frac{\Delta(t)dt}{t^2}.$$

Let us first not take into account the dependancies on $n$ and $g$.
As $$\rm Li\it(t^\rho)\sim \frac{t^\rho}{\log t^\rho},\ \ \ \rho<1, $$ $$ \int_2^x\rm Li\it(t^\rho)t^{-2}dt\ \ \ \text{is convergent}.$$

In order to prove the convergence of the second term, we need the following lemma:
\begin{lemme} \label{C2}
For all $x$ such that $$ (C2)\qquad\ \log x\geq 32^2 c_2^{-2}n\log^2 \frac{n^\frac{1}{2}}{c_2},$$ we have $$\exp \left(-c_2n^{-\frac{1}{2}}\log^{\frac{1}{2}}x\right)\leq \log^{-2}x.$$
\end{lemme}

\begin{preuvelemme} 
Put $y=\log^{\frac{1}{2}} x $ et $c=n^\frac{1}{2}/c_2$. Consider $f(y)=y^4 \exp -\frac{y}{c}$. We have to prove that $f(y)$ is less than $1$ if $y$ is big enough. We prove that $f$ is decreasing for $y\geq 4c$. Assume first that $c\leq e$. Then $y=16c$ verifies the inequality $f(y)\leq 1$. Indeed $$f(16c)=2^{16} c^4 e^{-16}=\frac{2^{16}}{e^{12}}\frac{c^4}{e^4}\leq 1.$$ 
If $c> e$, then $y=32c\log c$\ fits. Indeed, 
$$f(y)=32^4c^4\log^4(c)\, e^{-32\log c}=\frac{2^{20}}{c^{24}}\frac{\log^4 c}{c^4}\leq 1,$$
this estimate finishing the proof. 
\end{preuvelemme}

Therefore we have $\exp{\left(-c_2\,n^{-\frac{1}{2}}\log^\frac{1}{2}(t)\right)}=\mathcal{O}((\log{t})^{-2})$, and $$ \int_2^x\frac{t\exp{\left(-c_2\,n^{-\frac{1}{2}}\log^\frac{1}{2}(t)\right)}dt}{t^2}\ \ \ \text{is convergent}.$$ We obtain the convergence of the second integral with $x\to\infty$, because, for $t$ big enough:
$$|\Delta(t)|\leq \rm Li\it(t^\rho)+c_1 t \exp{\left(-c_2 n^{-{\frac{1}{2}}}\log^\frac{1}{2}(t)\right)}.$$

Finally we get $$\sum_{\mathcal{N}P\leq x} \frac{1}{\mathcal{N}P}=\log\log{x}+B+o(1).$$
\end{preuve}


We make now the residue appear continuing the asymptotical expansion of $\sum\frac{1}{\mathcal{N}P}$ by the calculation of the constant term $B.$

\begin{prop}
$$B=\sum_P\left\{\log\left(1-\frac{1}{\mathcal{N}P}\right)+\frac{1}{\mathcal{N}P}\right\}+\gamma+\log{\rm{Res}\it _{s=1}\zeta_K(s)}.$$
\end{prop}

\begin{preuve} For a complete proof, we refer to the article of Rosen \cite{ROS}, but let us still give the sketch of the proof. Write $C(x)=\log\log{x}+B+\epsilon(x).$
For $\delta>0$ define $$g(\delta)=\sum_{P} \frac{1}{\mathcal{N}P^{1+\delta}},$$and $$f(\delta)=g(\delta)-\log\zeta(1+\delta).$$
After some computation using the Abel transform, we find $g(\delta)=B-\gamma-\log{\delta}+\mathcal{O}(\delta)$.
Comparing with $\log{\zeta_K(1+\delta)}=-\log{\delta}+\log{\rm{Res}\it _{s=1}\zeta_K(s)}+\mathcal{O}(\delta),$ we get $f(\delta)=B-\gamma-\log{\rm{Res}\it _{s=1}\zeta_K(s)}+\mathcal{O}(\delta)$. Taking the limit when $\delta\to 0$, we obtain:
$$B=\sum_P\left\{\log\left(1-\frac{1}{\mathcal{N}P}\right)+\frac{1}{\mathcal{N}P}\right\}+\gamma+\log{\rm{Res}\it _{s=1}\zeta_K(s)}.$$ \end{preuve} 

We finally conclude that $$\sum_{\mathcal{N}P\leq x}\log{\frac{\mathcal{N}P}{\mathcal{N}P-1}}=\log\log{x}+\gamma+\log{\rm{Res}\it _{s=1}\zeta_K(s)}+o(1).$$ 

Let us now estimate the error term $$\epsilon(x)=\Delta(x)x^{-1}-\int_x^\infty\Delta(t)t^{-2}dt,$$ as the  function of $n$ and $g.$ 

\begin{prop} There are computable constants such that:
\begin{align*}(GRH)\qquad\text{For any }  x\geq 2\text{ we have } |\epsilon(x)|\leq & c\, x^{-\frac{1}{2}}\left(6g+3n\log{x}+2\,n\right),\\
 \text{For } x>>1,   \text{ we have }|\epsilon(x)|\leq &c_4 \frac{1}{\rho\log{x}}\left(1+(1-\rho)^{-1}\right)\\ &+2c_1\log^{-1}{x},
 \end{align*}
$x>>1$ meaning that $x$ must verify conditions $(C1)$and $(C2)$, the term in $\rho$ being present only if $\zeta_K$ has an exceptional zero.
\end{prop}

\begin{preuve} Assuming GRH, we obtain directly:
$$|\epsilon(x)|\leq c x^{-\frac{1}{2}}\left(2g+n\log{x}\right)+\int_x^\infty c\, t^{-\frac{3}{2}}\left(2g+n\log{t}\right)dt$$
$$|\epsilon(x)|\leq c x^{-\frac{1}{2}}\left(2g+n\log{x}\right)+2c x^{-\frac{1}{2}} (n\log{x}+2\,n+2g),$$
and finally:
$$(GRH)\qquad \qquad \ |\epsilon(x)|\leq c x^{-\frac{1}{2}}\left(6g+3n\log{x}+4\,n\right).$$

If we do not believe in GRH, we have to use the prime ideal theorem again : for $x$ verifying $(C1)$,
$$|\Delta(x)|\leq \rm Li\it(x^\rho)+c_1 x \exp{\left(-c_2\,n^{-\frac{1}{2}}\log^\frac{1}{2}(x)\right)},$$ 
 
Consider first the term $\Delta_1:=\rm Li\it(x^\rho).$ Put
$$\epsilon_1(x)=\Delta_1(x)x^{-1}-\int_x^\infty\Delta_1(t)t^{-2}dt.$$ 
Let $c_4$ be a constant such that $\rm Li\it(x)\leq c_4 x\log^{-1}{x}$ (for example $(1-\log2)^{-1}$). 
One has $$|\epsilon_1(x)|/c_4\leq \frac{1}{\rho\log{x}}+\int_x^\infty \frac{dt}{\rho\ t^{2-\rho}\log{t}}.$$ 
We can then easily bound the first error term by the following:

$$|\epsilon_1(x)|/c_4\leq \frac{1}{\rho\log{x}}\left(1+(1-\rho)^{-1}\right).$$

\paragraph{}

We now have to deal with the second error term $$\Delta_2(x)= c_1 x \exp{\left(-c_2\,n^{-\frac{1}{2}}\log^\frac{1}{2}(x)\right)}.$$

Using Lemma \ref{C2}, for $x$ verifying the condition $(C_2)$, we have:$$\exp{\left(-c_2\,n^{-\frac{1}{2}}\log^\frac{1}{2}{x}\right)}\leq \log^{-2}{x}.$$ Put
$$\epsilon_2(x)=\Delta_2(x)x^{-1}-\int_x^\infty\Delta_2(t)t^{-2}dt,$$
Thus, for $x$ verifying the  hypotheses $(C1)$ and $(C2)$ (note that condition $(C2)$ is very weak as compared to condition $(C1)$) and $x\geq e$, we obtain:
 $$\epsilon_2(x)\leq c_1(\log^{-2}{x}+\log^{-1}{x})\leq 2c_1\log^{-1}{x}.$$ \end{preuve}

\subsubsection*{End of the proof of the Mertens theorem:}

Let us start with the equality
$$\sum_{\mathcal{N}P\leq x} \frac{1}{\mathcal{N}P}=\log\log{x}+B_K+\epsilon_K(x).$$
One has
\begin{align*}
\sum_{\mathcal{N}P\leq x}\log\left(\frac{\mathcal{N}P}{\mathcal{N}P-1}\right)=&\log\log{x}+\sum_{\mathcal{N}P>x}\left\{\log\left(1-\frac{1}{\mathcal{N}P}\right)+\frac{1}{\mathcal{N}P}\right\}\\ 
&+\gamma+\log{\rm{Res}\it _{s=1}\zeta_K(s)}+\epsilon_K(x).
\end{align*}

We can bound the remainder term in the following way:$$\left|\sum_{\mathcal{N}P>x}\left\{\log\left(1-\frac{1}{\mathcal{N}P}\right)+\frac{1}{\mathcal{N}P}\right\}\right|\leq \sum_{\mathcal{N}P>x}\frac{1}{\mathcal{N}P^2}.$$

This sum can be calculated easily under GRH using the prime ideal theorem:

\begin{align*}
D(x)=\sum_{\mathcal{N}P> x} \frac{1}{\mathcal{N}P^2}=\int_x^\infty \frac{dt}{t^2\log{t}}+\int_x^\infty t^{-2}d\Delta(t)\leq&\frac{1}{x\log{x}}+\frac{|\Delta(x)|}{x^2}\\ &+2\int_x^\infty |\Delta(t)|t^{-3}dt,\end{align*}
therefore
$$(GRH)\qquad  \text{for any } x\geq 2 \text{ we have } D(x)\leq\frac{1}{x\log{x}}+\frac{10\,g+3n\log x}{3 \,x\sqrt{x}}+\frac{2n}{x}.$$

Without GRH,  we can use the bound for $\pi(x)$ (see \cite{SC}) valid for $$ (C3) \qquad\qquad \log{x}\geq c_5\,g\log{2g}\log\log{12g}:$$
$$\pi(x)\leq c_6 x\log^{-1}(x).$$

We have $$D(x)=\sum_{\mathcal{N}P> x} \frac{1}{\mathcal{N}P^2}=\int_{x}^\infty\frac{d \pi(t)}{t^2}$$ and, for $x$ sufficiently large,
$$D(x)=-\frac{\pi(x)}{x^2}+2\int_x^\infty\pi(t)t^{-3}dt\leq 2\,c_6\int_x^\infty t^{-2}\log^{-1}(t)dt\leq \frac{2\,c_6}{x\log{x}}.$$

Putting all this together, we obtain the following. For $x$ verifying $(C_1)$, $(C_2)$ and $(C_3)$
\begin{align*}\sum_{\mathcal{N}P\leq x}\log\left(\frac{\mathcal{N}P}{\mathcal{N}P-1}\right)=&\log\log{x}+\gamma+\log{\rm{Res}\it _{s=1}\zeta_K(s)}\\ &+\mathcal{O}\left(\frac{1}{\log{x}}\right)+\frac{1}{1-\rho}\mathcal{O}\left(\frac{1}{\log{x}}\right), \end{align*} 
where the term in $\rho$ is there only if  $K$ has an exceptional zero. The classical Mertens theorem follows by an easy application of the Taylor expansion.

Under GRH, we obtain a stronger result true for $x\geq 2$:
\begin{align*}\sum_{\mathcal{N}P\leq x}\log\left(\frac{\mathcal{N}P}{\mathcal{N}P-1}\right)=&\log\log{x}+\gamma+\log{\rm{Res}\it _{s=1}\zeta_K(s)}\\ & +n\,\mathcal{O}\left(\frac{\log{x}}{\sqrt{x}}\right)+g\,\mathcal{O}\left(\frac{1}{\sqrt{x}}\right),\end{align*}
 which leads to the Mertens theorem under GRH.

\subsection{Proof in the case of algebraic varieties}

We first establish the Mertens theorem in the case of smooth absolutely irreducible projective algebraic varieties. The generalised Brauer--Siegel follows immediately from it.

For any sequence $(v_n)$ such that the radius of convergence $\rho$ of the series $\sum v_n t^n$ is strictly positive, put $$\psi_{m,v}(t)=\sum_{n=1}^{+\infty}v_{mn} t^{mn},$$ and $\psi_v(t)=\psi_{1,v}(t)$. For $t<r^{-d}\rho$, we have the explicit formulae:
\begin{theo}[Explicit Formula]
$$\sum_{f=1}^{+\infty}f\Phi_{r^f}\psi_{f,v}=\psi_v(t)+\psi_v(r^d t)+\sum_{i=1}^{2d-1}(-1)^i\sum_{j=1}^{b_i}\psi_v(r^{\frac{i}{2}}\omega_{i,j}t)$$
\end{theo}
\begin{preuve}\cite{TL} \end{preuve}	

Choose $N\in\mathbb{N}$ et take $v_n(N)=\frac{1}{n}$ if $n\leq N$ and $0$ otherwise. Applying this explicit formula with $t=r^{-d}$, we get : $$S_0(N)=S_1(N)+S_2(N)+S_3(N),$$ where
\begin{align*}
S_0(N)= &\sum_{n=1}^{N}n^{-1} r^{-dn} \sum_{m/n}m\Phi_{r^{dm}},\\
S_1(N)= & \sum_{n=1}^N\frac{1}{n},\\
S_2(N)= & \sum_{n=1}^N\frac{1}{nr^{dn}},\\
S_3(N)= & \sum_{i=1}^{2d-1}(-1)^i\sum_{j=1}^{b_i}\sum_{n=1}^N\frac{1}{n}(r^{\frac{i}{2}-d}\omega_{i,j})^n.
 \end{align*}


\begin{lemme} $$0\leq \sum_{f=1}^N \Phi_{r^f} \log\frac{r^{df}}{r^{df}-1}-S_0(N)\leq  \frac{8}{N\,r^{dN/2}}+\frac{6\,b}{N\,r^{(d+\frac{1}{2})\frac{N}{2}}}.$$
\end{lemme}
\begin{preuvelemme} Let us first transform the expression of $S_0$:
\begin{align*}
S_0(N)= & \sum_{f=1}^{N}\sum_{m=1}^{E(N/f)}f\Phi_{r^f} r^{-dfm}(fm)^{-1}\\
= & \sum_{f=1}^{N}\Phi_{r^f} \sum_{m=1}^{E(N/f)}\frac{1}{r^{dfm}m}.
\end{align*}
Then evaluate $S_0$:
\begin{align*}
0\leq \sum_{f=1}^N \Phi_{r^f} \log\frac{r^{df}}{r^{df}-1}-S_0(N)= & \sum_{f=1}^{N}\Phi_{r^f} \left(\log\frac{r^{df}}{r^{df}-1}-\sum_{m=1}^{E(N/f)}\frac{1}{r^{dfm}m}\right) \\
= & \sum_{f=1}^{N}\Phi_{r^f}\sum_{m=E(N/f)+1}^\infty\frac{1}{r^{dfm}m}.
\end{align*}
As $1/m\leq 1/(E(N/f)+1),$ we get:

$$0\leq \sum_{f=1}^N \Phi_{r^f} \log\frac{r^{df}}{r^{df}-1}-S_0(N)\leq \sum_{f=1}^{N}\frac{\Phi_{r^f}}{(E(N/f)+1)(r^{df})^{E(N/f)}(r^{df}-1)}.$$

In order to deal with $\Phi_{r^f}$ we use $$\Phi_{r^f}\leq \frac{r^{df}+1+\sum_{i=1}^{2d-1}r^{if/2}b_i}{f}.$$ Let $b=b_X=\max_i(b_i).$ We obtain:

\begin{eqnarray*}
 \Phi_{r^f} & \leq & \frac{1}{f}\left(r^{df}+1+b\sum_{i=1}^{2d-1} r^{if/2} \right) \\
& \leq & \frac{1}{f}\left(r^{df}+1+b\, r^{\frac{f}{2}} \frac{r^{\frac{2d-1}{2}f}-1}{r^{\frac{f}{2}}-1}\right) \\
& \leq & \frac{1}{f}\left(r^{df}+1+2\,b\, r^{df-\frac{f}{2}}\right).
\end{eqnarray*}
Thus 
$$0\leq\sum_{f=1}^N \Phi_{r^f} \log\frac{r^{df}}{r^{df}-1}-S_0(N)\leq \frac{1}{N}\sum_{f=1}^{N}\frac{\left(r^{df}+1+2 b\, r^{df-\frac{f}{2}}\right)\left(r^{df}-1\right)^{-1}}{r^{dfE(N/f)}}.$$ We split our sum in two in the following way: for $f>E(N/2)$ where $E(N/f)=1$, and for $f\leq E(N/2)$ where we use $fE(N/f)\leq N-f.$ 

\begin{eqnarray*}
{0\leq\sum_{f=1}^N \Phi_{r^f} \log(\frac{r^{df}}{r^{df}-1}) -S_0(N)}& \leq &\frac{1}{N}\sum_{f=1}^{E(N/2)}\frac{2+4\,b\,r^{-\frac{f}{2}}}{r^{d(N-f)}}\\ & & +\frac{1}{N}\sum_{f>E(N/2)}^N \frac{2+ 4\,b\,r^{-\frac{f}{2}}}{r^{df}}\\
& \leq & \frac{8+12\,b\,r^{-\frac{N}{4}}}{N\,r^{dN/2}}.
\end{eqnarray*}

We finally obtain the following inequality:
$$0\leq \sum_{f=1}^N \Phi_{r^f} \log\frac{r^{df}}{r^{df}-1}-S_0(N)\leq  \frac{8}{N\,r^{dN/2}}+\frac{6\,b}{N\,r^{(d+\frac{1}{2})\frac{N}{2}}}.$$
\end{preuvelemme}

In order to estimate $S_1$ we use the following well-known inequality (see \cite{HW}):
\begin{lemme} $$\frac{1}{N(N+1)}\leq S_1(N)-\log N-\gamma\leq\frac{1}{N}.$$
\end{lemme}


\begin{lemme} $$0\leq \log\frac{r^d}{r^d-1}-S_2(N)=\sum_{n=N+1}^\infty\frac{1}{nr^{dn}}\leq \frac{1}{r^{dN}(N+1)(r^d-1)}.$$
\end{lemme}
\begin{preuvelemme} $S_2$ is the partial summation of the entire function $\log\frac{r^d}{r^d-1}.$ The inequality comes from the estimation of the remainder term.
\end{preuvelemme}

Let us recall first that:$$\log{\rm{Res}\it _{s=d}(\log{r}\,\zeta(s))}-\log\frac{r^d}{r^d-1}=\sum_{i=1}^{2d-1}(-1)^{i+1}\sum_{j=1}^{b_i}\log{\left(1-r^{\frac{i}{2}-d}\omega_{i,j}\right)}.$$ Compute now $S_3:$

\begin{lemme} $$|S_3(N)-\sum_{i=1}^{2d-1}(-1)^{i+1}\sum_{j=1}^{b_i}\log{\left(1-r^{\frac{i}{2}-d}\omega_{i,j}\right)}|\leq\frac{b}{(r^{\frac{1}{2}}-1)(N+1)(r^\frac{N}{2}-1)}.$$
\end{lemme}
\begin{preuvelemme} Consider $$R(N)=|S_3(N)-\sum_{i=1}^{2d-1}(-1)^{i+1}\sum_{j=1}^{b_i}\log{\left(1-r^{\frac{i}{2}-d}\omega_{i,j}\right)}|.$$ One has
$$R(N)=|\sum_{i=1}^{2d-1}(-1)^i\sum_{j=1}^{b_i}\sum_{n=N+1}^\infty\frac{1}{n}(r^{\frac{i}{2}-d}\omega_{i,j})^n|,$$ therefore

\begin{align*}
R(N) & \leq
\sum_{i=1}^{2d-1}\sum_{j=1}^{b_i}\frac{1}{N+1}\sum_{n=N+1}^\infty(r^{\frac{i}{2}-d})^n,\\
& \leq\frac{1}{N+1}\sum_{i=1}^{2d-1} \frac{b_i(r^{\frac{i}{2}-d})^{N+1}}{1-r^{\frac{i}{2}-d}},
\end{align*}
and
$$R(N)\leq\frac{b}{(r^{\frac{1}{2}}-1)(N+1)r^{dN}}\sum_{i=1}^{2d-1}r^{\frac{iN}{2}},$$ which leads to the result.
\end{preuvelemme}

\subsection*{Putting everything together:}

We deduce then, that for $N$ big enough,  \begin{align*}
\log \prod_{f=1}^{N}\left(1-\frac{1}{r^{df}}\right)^{\Phi_{r^f}}=& -\log N-\gamma+\log\left(1-\frac{1}{r^d}\right)-\log\left(1-\frac{1}{r^d}\right)\\ 
& -\log{(\log{r} \rm{Res}\it _{s=d}\zeta_X)} +\mathcal{O}\left(\frac{1}{N}\right)+b\,\mathcal{O}\left(\frac{r^\frac{-N}{2}}{N}\right),
\end{align*} 
where the constant involved in $\mathcal{O}$ does not depend on $K$.  \hfill $\square$

\section{Proof of the generalised BS theorem}

\subsection{Without GRH}

If we do not believe in the generalised Riemann hypothesis, we have to take into account the conditions on $x$ which forbid us to take the limit. Consider an asymptotically exact family $(K_i)$ of number fields, and divide it into three subfamilies. The first one consists in all the fields that have no exceptional zeroes, we include in the second one the fields that do have exceptional zeroes but no Siegel zero, and the last one contains the fields that have a Siegel zero. If one of them is finite, we omit it. 

Let us focus on the second and the third families, the case of the first one being much easier because of the absence of the $\rho$-term (or take $\rho=0$ in the following).
Let us specialise the Mertens theorem in $x=e^{C\,n\,g^2(1-\rho)^{-1}}$, where $C$ is big enough to allow $x$ verify all the three conditions. Thus, for $g$ big enough and $M$ an explicit constant:
$$\left|\sum_{q\leq e^{C\,n\,g^2(1-\rho)^{-1}}}\frac{\Phi_q}{g}\log\left(\frac{q}{q-1}\right)-\frac{\log{\varkappa}}{g}\right|\leq M \frac{\log{g}}{g}-\frac{\log(1-\rho)}{g}.$$
\begin{lemme} Consider the family $(K_i)$ and its exceptional zeroes $\rho_i.$ Suppose that $$\lim_i\log(1-\rho_i)/g_i=0.$$ Then $\kappa$ exists and we have $$\kappa=\sum_q\phi_q\log\frac{q}{q-1}.$$
\end{lemme}
Let us assume the lemma. Look first at the second subfamily still denoted by $(K_i).$ Each  $\zeta_{K_i}$ has an exceptional zero verifying $1-(8g)^{-1}\leq \rho<1-(32g)^{-1},$ thus $(1-\rho)^{-1}\leq 32g.$ Taking the logarithm, we see that this family verifies the condition of the lemma.

The case of the third subfamily, which is still denoted by $(K_i)$ for the sake of commodity, is not so easy, because $\rho$ can go very close to $1.$ In order to control the magnitude of the term in $\rho$, we need to assume that the fields are almost normal (or some additional condition as below). Indeed, thanks to Stark we know  that a Siegel zero $\rho$ of an almost normal number field $K$ is also a Siegel zero of a subextension of  $K$ of degree $2$ over $\mathbb{Q}$ (see \cite{Sta}). In addition, we can estimate $(1-\rho)^{-1}$ as follows \cite{LOU}:

\begin{lemme} Let $K$ be a number field of degree $n_K>1$. Then $$\frac{1}{1-\rho_K}\leq \varkappa_K^{-1}\left(\frac{g_K}{n_K}\right)^{n_K}.$$
\end{lemme}

Let $(k_i)$ be a family of quadratic extensions of $\mathbb{Q}$ having the same Siegel zeroes as $(K_i)$ and let us apply this lemma :
$$-\log (1-\rho_{K_i})=-\log(1-\rho_{k_i})\leq -\log \varkappa_{k_i}+2\log\left(\frac{g_{k_i}}{2}\right).$$

Thus we obtain $$0<-\frac{\log (1-\rho_{K_i})}{g_{K_i}}\leq -\frac{\log \varkappa_{k_i}}{g_{K_i}}+2g_{K_i}^{-1}\log\frac{g_{k_i}}{2}.$$

As $k_i\subset K_i$, $g_{K_i}\geq g_{k_i}$; the first and the last term of the right side of the inequality tend to zero with $i.$  The second term, if positive, can be bounded by $g_{k_i}^{-1}\log\kappa_i$ and we use the classical Brauer--Siegel theorem for quadratic fields which says that it tends to $0$. We then apply the lemma to deduce the generalised Brauer--Siegel theorem. 
\paragraph{}

We still have to prove the first lemma.

\begin{preuvelemme}
Put $$f_g(q)=\frac{\Phi_q}{g}\log\left(\frac{q}{q-1}\right)\delta_g(q),$$ where $\delta_g(q)=1$ if $q\leq e^{C\,n\,g^2(1-\rho)^{-1}}$, $0$ otherwise. 
Now let the genus tend to infinity ($g_{K_i}$ being $g_i$ again).
As $$\sup_{i>>1}\sum_q f_{g_i}(q)\leq\sup_{i}\frac{\log{\varkappa_{K_i}}}{g_i}+1,$$
this last quantity being well defined because of basic inequality of \cite{TVF}, and we can apply the Fatou lemma, and obtain:
\begin{align*}
\sum_q \phi_q \log\left(\frac{q}{q-1}\right)&=\sum_q\liminf_{i\to\infty}f_{g_i}(q)\\
&\leq\liminf_{i\to\infty}\sum_q{f_{g_i}(q)}=\liminf_{i\to\infty}\frac{\log{\varkappa_{K_i}}}{g_i}.\end{align*}
Combining this result with the Brauer--Siegel inequality (theorem \ref{IBS}), we deduce the existence of the limit of  
$\frac{\log{\varkappa_{K_i}}}{g_i}$, which equals $\sum_q \phi_q \log\left(\frac{q}{q-1}\right)$
\end{preuvelemme}

Remark that in the case of the first subfamily, the proof becomes easier, because we do not have to deal with the $\rho$-term. The bound in the Mertens theorem is only in $\log{g}/g.$ Specialising in $x=e^{C\,n\,g^2}$ instead, and suppressing the $\rho$-term in $f_g(q)$ ($q>e^{C\,n\,g^2}$), we obtain the generalised Brauer--Siegel theorem. \hfill \mbox{$\Box$}

Let us now prove Proposition \ref{propcomp}:\\ i. Recall the following key-lemma of Stark.
\begin{lemme}[\cite{Sta} Lemma 8] Let $k$ be a number field of degree $n_k>1.$ Assume that there is a $\beta\in \mathbb{R}$ such that $$1-\frac{1}{8n_k!g_k}\leq\beta<1$$ and $\zeta_k(\beta)=0.$ Then there is a quadratic subfield $F$ of $k$ such that $\zeta_F(\beta)=0.$
\end{lemme}
Assume as before that $(K_i)$ has an infinite number of Siegel zeroes which do not verify the condition of the lemma.  

Let us split as before the family $(K_i)$ into three subfamilies. The first one containing the fields that do not have an exceptional zero, the second one consisting in the fields that have zeroes that do not verify this lemma and the last one consisting in those whose Siegel zeroes verify the condition of the lemma. If one of these families is finite, we omit it. The first and the third cases have already been treated before, so let us consider the second subfamily. Let us call it $(K_i)$ again. We still have to bound $g^{-1}\log\frac{1}{1-\rho}.$ Their exceptional zeroes verify $$\log\frac{1}{1-\rho}\leq 8+\log{n!}+\log{g}.$$ As $\log{n!}\leq n\log{n}$, we deduce that $$\frac{1}{g}\log\frac{1}{1-\rho}\leq \frac{n\log{n}}{g}+m\frac{\log{g}}{g},$$ where $m$ is an explicit positive constant. Therefore this quantity tends to $0$ and this completes the proof.

ii. Suppose now that the family $(K_i)$ is asymptotically good, and that an infinite number of them admit a Siegel zero. Then, because of Louboutin's lemma, we obtain $$\frac{1}{g_i}\log\frac{1}{1-\rho}\leq -\frac{\log\varkappa_i}{g_i}+\frac{n_i}{g_i}\log{e\frac{g_i}{n_i}}.$$ This leads to the result, since $\phi_\infty=\lim\frac{n_i}{g_i}.$ \hfill \mbox{$\Box$}

\subsection{Assuming GRH}
In the number field case, let $(K_i)$ be a family of fields with $g_i\to\infty.$ Starting with the Mertens theorem, dividing by $g_i$, taking $g_i\to\infty$ (we can do it because this time there is no condition on $x$), we obtain the Brauer--Siegel theorem. Indeed, the last paragraph shows that the limit of $\log\varkappa_{K_i}/g_i$ exists, and the asymptotical result follow directly.

In the variety case, let $(X_i)$ be a family of smooth absolutely irreducible projective algebraic varieties over $\mathbb{F}_r$. We assume either the result of Zykin, or use the bounds for $\Phi_{r^f}$ that we needed for the Mertens theorem, in order to prove that the series $$\sum \phi_{r^m}\log\frac{r^{dm}}{r^{dm}-1}\ \text{ is convergent,}$$ and that
$$\lim_{b\to\infty}\frac{1}{b}\sum_{m=1}^{f(b)}\Phi_{r^m}\log\frac{r^{dm}}{r^{dm}-1}=\sum_{m=1}^\infty\phi_{r^m}\log\frac{r^{dm}}{r^{dm}-1}, $$ for any function $f$ of $b$ verifying the conditions :
$$\begin{cases}\displaystyle \lim_{b\to\infty} f(b)=\infty,\\
\displaystyle\lim_{b\to\infty}\frac{f(b)}{b}=0.\end{cases}$$ This is a bit technical but not hard (for the case $d=1$ see \cite{TIF}).

Using this result in the Mertens theorem (we put $N=f(b),$ divide by $b$ and make $b\to\infty$) gives us that the limit of $\log\varkappa_{X_i}/b_{X_i}$ exists.  We divide now by $b_{X_i}$ (for any $N$) in the Mertens theorem and make $b_{X_i}\to\infty$ in order to obtain our version of the Brauer--Siegel theorem for varieties.

One could likely obtain similar results in the non-smooth case, using the virtual Betti numbers.  We hope to do this in further work. Let us conclude by the following remark. The explicit Mertens theorem is much more interesting than its application to the generalised Brauer--Siegel theorem, because it contains more information, and can be therefore useful, for example if we would like to look at the problem in the classical way, putting all our attention to the residues $\varkappa_i$ instead of to the convergent series, and consider the limit of $\varkappa_i/g_i$ in the tower.

\section*{Acknowledgements}
I would like to thank Gilles Lachaud, Michael Tsfasman and Alexei Zykin for very useful discussions, Mireille Car for letting me know her results on the Mertens  theorem and Michel Balazard for having pointed out some mistakes in the first version of this paper.

\bibliographystyle{amsplain}
\bibliography{biblio}

\end{document}